\documentclass{amsproc}%
\usepackage{graphicx}
\usepackage{amscd}
\usepackage{amsmath}
\usepackage{amsfonts}
\usepackage{amssymb}%
\theoremstyle{plain}

\numberwithin{equation}{section}

\begin{document}
\title[A half-space approach to order dimension]{A half-space approach to order dimension}
\author{Stephan Foldes}
\address{Institute of Mathematics, Tampere University of Technology PL 553, 33101
Tampere, Finland}
\email{stephan.foldes@tut.fi}
\author{Jen\H{o} Szigeti}
\address{Institute of Mathematics, University of Miskolc, Miskolc, Hungary 3515}
\email{jeno.szigeti@uni-miskolc.hu}
\thanks{\noindent The work of the first named author was partially supported by the
European Community's Marie Curie Program (contract MTKD-CT-2004-003006). The
second named author's work was supported by OTKA of Hungary No. T043034.}

\begin{abstract}
The aim of the present paper is to investigate the half-spaces in the
convexity structure of all quasiorders on a given set and to use them in an
alternative approach to classical order dimension. The main result states that
linear orders can almost always be replaced by half-space quasiorders in the
definition of the dimension of a partially ordered set.

\end{abstract}
\maketitle

\noindent1. INTRODUCTION

\bigskip

\noindent Within the framework of the general theory of abstract convexity
(van de Vel [9]), strict quasiorders (irreflexive and transitive relations) on
a set $A$ can be thought of as convex subsets of $\{(x,y)\in A\times A\mid
x\neq y\}$:

\begin{enumerate}
\item $\{(x,y)\in A\times A\mid x\neq y\}$ is a strict quasiorder,

\item any intersection of strict quasiorders is a strict quasiorder,

\item any nested union of strict quasiorders is a strict quasiorder.
\end{enumerate}

\noindent In general, a half-space is defined as a convex subset of the base
set with a convex set complement. Abstract convexity theory addresses
questions such as the representation of convex sets as intersections of
half-spaces. For technical reasons, instead of the strict quasiorders in
$\{(x,y)\in A\times A\mid x\neq y\}$, we shall consider the ordinary
(reflexive) quasiorders in $A\times A$ (there is a natural one to one
correspondence between them). We can use half-space quasiorders to define the
half-space dimension of a quasiordered set, in a similar way as linear orders
are used to define the order dimension of a partially ordered set. The aim of
the present paper is to investigate the half-space quasiorders and to study
the above dimension concept for quasiorders, along the lines of the classical
theory of order dimension (see e.g. [1,2,7,8]). Our main result (Theorem 2.16)
states that linear orders can almost always be replaced by half-space
quasiorders in the definition of the order dimension. Since there are
considerably more half-spaces than linear orders, establishing upper bounds on
order dimension can be easier using representations of partial orders as
intersections of half-spaces.

\noindent In section 2 we provide some simple characterizations of half-spaces
and examine the relationship between half-spaces and linear orders. A standard
construction together with a complete description of half-spaces is also
given. In the rest of section 2, we show the tight connection between
half-space dimension and classical order dimension. It turns out, that the
half-space dimension and the order dimension of a partially ordered set can be
different only for half-space partial orders.

\noindent In section 3 we prove that the direct product of quasiorders can be
a half-space only in one exceptional situation.

\bigskip

\noindent2. HALF-SPACES\ AND\ THE DIMENSION\ OF\ QUASIORDERED\ SETS

\bigskip

\noindent A \textit{quasiorder} $\gamma$ on the set $A$ is a reflexive and
transitive relation:
\[
\Delta_{A}=\{(a,a)\mid a\in A\}\subseteq\gamma\subseteq A\times A
\]
and $(x,y)\in\gamma$, $(y,z)\in\gamma$ imply $(x,z)\in\gamma$ for all
$x,y,z\in A$. The containment relation $\subseteq$ provides a natural complete
lattice structure on the set Quord$(A)$\ of all quasiorders on $A$:
$($Quord$(A),\vee,\cap)$. If $\gamma$ is a partial order, then we frequently
use the standard notations $x\leq_{\gamma}y$ and $x<_{\gamma}y$ for
$(x,y)\in\gamma$ and for $(x,y)\in\gamma$, $x\neq y$. For a quasiorder
$\gamma$, the relation $\gamma\cap\gamma^{-1}$ is an equivalence on $A$, the
equivalence class of an element $a\in A$ is denoted by $[a]_{\gamma\cap
\gamma^{-1}}$, thus
\[
A/(\gamma\cap\gamma^{-1})=\{[a]_{\gamma\cap\gamma^{-1}}\mid a\in A\}.
\]
It is well known that $\gamma$ induces a natural partial order $r_{\gamma}$
(in order to avoid repeated indices, we write $\leq^{\gamma}$ instead of
$\leq_{r_{\gamma}}$) on the above quotient set: for $a,b\in A$%
\[
\lbrack a]_{\gamma\cap\gamma^{-1}}\leq^{\gamma}[b]_{\gamma\cap\gamma^{-1}%
}\text{ if and only if }(x,y)\in\gamma\text{ for some }x\in\lbrack
a]_{\gamma\cap\gamma^{-1}}\text{ and }y\in\lbrack b]_{\gamma\cap\gamma^{-1}%
}\text{.}%
\]
Also $[a]_{\gamma\cap\gamma^{-1}}\leq^{\gamma}[b]_{\gamma\cap\gamma^{-1}}$
holds if and only if $(x,y)\in\gamma$ for all $x\in\lbrack a]_{\gamma
\cap\gamma^{-1}}$ and for all $y\in\lbrack b]_{\gamma\cap\gamma^{-1}}$.

\noindent A quasiorder $\alpha\subseteq A\times A$ is said to be a
\textit{half-space on }$A$ if it has a ''strong'' complement in the lattice
$($Quord$(A),\subseteq)$, i.e. if $\alpha\cup\beta=A\times A$ and $\alpha
\cap\beta=\Delta_{A}$ hold for some quasiorder $\beta\subseteq A\times A$.
Clearly, this complement $\beta$ is also a half-space and is uniquely
determined by $\alpha$: $\beta=\Delta_{A}\cup((A\times A)\setminus\alpha)$. It
follows, that $\alpha$ is a half-space if and only if $\Delta_{A}\cup((A\times
A)\setminus\alpha)$ is transitive. The simplest examples of half-spaces are
linear orders, the identity $\Delta_{A}$\ and the full relation $A\times A$ on
any set $A$. Complementary half-spaces are put into a pair of the form
$\alpha\updownarrow\beta$ and can be characterized in the lattice
$($Quord$(A),\vee,\cap)$ as follows.

\bigskip

\noindent\textbf{2.1.Proposition.}\textit{\ For any quasiorders }$\alpha
,\beta\in$Quord$(A)$\textit{\ the following are equivalent:}

\begin{enumerate}
\item $\alpha\updownarrow\beta$\textit{\ is a pair of complementary
half-spaces, i.e. }$\alpha\cap\beta=\Delta_{A}$\textit{\ and }$\alpha\cup
\beta=A\times A$\textit{.}

\item $\alpha\cap\beta=\Delta_{A}$\textit{\ and }$(\alpha\cap\gamma)\vee
(\beta\cap\gamma)=\gamma$\textit{\ for all }$\gamma\in$Quord$(A)$\textit{.}
\end{enumerate}

\bigskip

\noindent\textbf{Proof.} $(1)\Longrightarrow(2)$:
\[
\gamma=(A\times A)\cap\gamma=(\alpha\cup\beta)\cap\gamma=(\alpha\cap
\gamma)\cup(\beta\cap\gamma)\subseteq(\alpha\cap\gamma)\vee(\beta\cap
\gamma)\subseteq\gamma.
\]
$(2)\Longrightarrow(1)$: Suppose that $\alpha\cup\beta\neq A\times A$, then
$(a,b)\notin\alpha\cup\beta$ for some $a,b\in A$. Since $\gamma(a,b)=\Delta
_{A}\cup\{(a,b)\}$ is a quasiorder on $A$, we have
\[
(\alpha\cap\gamma(a,b))\vee(\beta\cap\gamma(a,b))=\gamma(a,b)
\]
in contradiction with $\alpha\cap\gamma(a,b))=\beta\cap\gamma(a,b)=\Delta_{A}
$.$\square$

\bigskip

\noindent For a half-space $\alpha$ the inverse relation $\alpha^{-1}$ is also
a half-space, if $\alpha\updownarrow\beta$ for $\alpha,\beta\in$Quord$(A)$,
then $\alpha^{-1}\updownarrow\beta^{-1}$. If $B\subseteq A$ is a subset, then
the restriction of a quasiorder to $B$ yields a quasiorder on $B$ and a
similar statement holds for half-spaces, $\alpha\updownarrow\beta$ implies
that $\alpha\cap(B\times B)\updownarrow\beta\cap(B\times B)$. This observation
can be used to give another characterization of half-spaces.

\bigskip

\noindent\textbf{2.2.Proposition.}\textit{\ For a quasiorder }$\alpha\in
$Quord$(A)$\textit{\ the following are equivalent:}

\begin{enumerate}
\item $\alpha$\textit{\ is a half-space.}

\item $\alpha\cap(B\times B)$\textit{\ is a half-space (on }$B$\textit{) for
any three element subset }$B\subseteq A$\textit{.}

\item \textit{For any }$x,y,z\in A$\textit{\ the relations }$(x,y)\notin
\alpha$\textit{, }$(y,x)\notin\alpha$\textit{\ and }$(x,z)\in\alpha$\textit{,
}$z\neq x$\textit{\ imply that }$(y,z)\in\alpha$\textit{.}

\item \textit{For any }$x,y,z\in A$\textit{\ the relations }$(z,y)\notin
\alpha$\textit{, }$(y,z)\notin\alpha$\textit{\ and }$(x,z)\in\alpha$\textit{,
}$x\neq z$\textit{\ imply that }$(x,y)\in\alpha$\textit{.}
\end{enumerate}

\bigskip

\noindent\textbf{Proof.} $(1)\Longrightarrow(2)$: This is a special case of
our claim preceding Proposition 2.2.

\noindent$(2)\Longrightarrow(3)$: Let $(x,y)\notin\alpha$, $(y,x)\notin\alpha
$\ and $(x,z)\in\alpha$, $z\neq x$ for the elements $x,y,z\in A$ and take the
three element subset $B=\{x,y,z\}$ of $A$. Suppose that $(y,z)\notin\alpha$
and consider the complementary half-space $\delta\subseteq B\times B$ of
$\alpha\cap(B\times B)$. Now
\[
(\alpha\cap(B\times B))\cup\delta=B\times B
\]
implies that $(x,y)\in\delta$ and $(y,z)\in\delta$, whence $(x,z)\in
(\alpha\cap(B\times B))\cap\delta=\Delta_{B}$ can be derived in contradiction
with $z\neq x$.

\noindent$(3)\Longrightarrow(4)$: Let $(z,y)\notin\alpha$, $(y,z)\notin\alpha$
and $(x,z)\in\alpha$, $x\neq z$ for the elements $x,y,z\in A$ and suppose that
$(x,y)\notin\alpha$. Clearly, $(y,x)\in\alpha$ would imply $(y,z)\in\alpha$, a
contradiction. Thus $(x,y)\notin\alpha$, $(y,x)\notin\alpha$ and
$(x,z)\in\alpha$, $x\neq z$, whence we obtain that $(y,z)\in\alpha$, a
contradiction. It follows that $(x,y)\in\alpha$.

\noindent$(4)\Longrightarrow(1)$: In order to see the transitivity of
$\beta=\Delta_{A}\cup((A\times A)\setminus\alpha)$ let $(x,y)\in\beta$,
$(y,z)\in\beta$, $x\neq y$ and suppose that $(x,z)\notin\beta$. We have either
$(z,y)\notin\alpha$ or $(z,y)\in\alpha$. In the first case $(z,y)\notin\alpha
$, $(y,z)\notin\alpha$ and $(x,z)\in\alpha$, $x\neq z$ would imply that
$(x,y)\in\alpha\cap\beta=\Delta_{A}$, a contradiction. In the second case
$(x,z)\in\alpha$ and $(z,y)\in\alpha$ would imply that $(x,y)\in\alpha
\cap\beta=\Delta_{A}$, a contradiction again. Thus we have $(x,z)\in\beta
$.$\square$

\bigskip

\noindent\textbf{2.3.Proposition.}\textit{\ If }$\alpha$\textit{\ is a
half-space quasiorder on }$A$\textit{, then the induced partial order
}$r_{\alpha}$\textit{\ is a half-space on }$A/(\alpha\cap\alpha^{-1}%
)$\textit{.}

\bigskip

\noindent\textbf{Proof.} We can use part (3) in Proposition 2.2. If
$([x]_{\alpha\cap\alpha^{-1}},[y]_{\alpha\cap\alpha^{-1}})\notin r_{\alpha}$,
$([y]_{\alpha\cap\alpha^{-1}},[x]_{\alpha\cap\alpha^{-1}})\notin r_{\alpha}%
$\ and $([x]_{\alpha\cap\alpha^{-1}},[z]_{\alpha\cap\alpha^{-1}})\in
r_{\alpha}$, $[z]_{\alpha\cap\alpha^{-1}}\neq\lbrack x]_{\alpha\cap\alpha
^{-1}}$, then we have $(x,y)\notin\alpha$, $(y,x)\notin\alpha$\ and
$(x,z)\in\alpha$, $z\neq x$. Since $\alpha$ is a half-space, we obtain first
$(y,z)\in\alpha$ and then $([y]_{\alpha\cap\alpha^{-1}},[z]_{\alpha\cap
\alpha^{-1}})\in r_{\alpha}$.$\square$

\bigskip

\noindent\textbf{2.4.Proposition.}\textit{\ If }$\gamma\subseteq A\times
A$\textit{\ is a quasiorder and }$\gamma\subseteq\alpha$\textit{\ for some
half-space }$\alpha$\textit{\ on }$A$\textit{, then there exists a half-space
}$\tau$\textit{\ on }$A$\textit{,\ such that }$\gamma\subseteq\tau
\subseteq\alpha$\textit{\ and }$\tau\cap\tau^{-1}=\gamma\cap\gamma^{-1}%
$\textit{.}

\bigskip

\noindent\textbf{Proof.} Let $R$ be a linear extension of the induced partial
order $r_{\gamma}$ and define the relation $\tau\subseteq A\times A$ as
follows:
\[
\tau=\alpha\setminus\{(a,b)\in\alpha\cap\alpha^{-1}\mid\lbrack b]_{\gamma
\cap\gamma^{-1}}<_{R}[a]_{\gamma\cap\gamma^{-1}}\}.
\]
Since $(x,y)\in\gamma$ implies that $(x,y)\in\alpha$ and $[x]_{\gamma
\cap\gamma^{-1}}\leq_{R}[y]_{\gamma\cap\gamma^{-1}}$, we obtain that
$(x,y)\in\tau$. Thus $\gamma\subseteq\tau\subseteq\alpha$ and $\gamma
\cap\gamma^{-1}\subseteq\tau\cap\tau^{-1}$. If $(x,y)\in\tau\cap\tau^{-1}$,
then the relations $[y]_{\gamma\cap\gamma^{-1}}<_{R}[x]_{\gamma\cap\gamma
^{-1}}$ and $[x]_{\gamma\cap\gamma^{-1}}<_{R}[y]_{\gamma\cap\gamma^{-1}}$ are
not satisfied, whence $[x]_{\gamma\cap\gamma^{-1}}=[y]_{\gamma\cap\gamma^{-1}%
}$ and $(x,y)\in\gamma\cap\gamma^{-1}$ can be derived. It follows, that
$\tau\cap\tau^{-1}\subseteq\gamma\cap\gamma^{-1}$ and hence $\tau\cap\tau
^{-1}=\gamma\cap\gamma^{-1}$.

\noindent In order to see the transitivity of $\tau$ take $(x,y)\in\tau$ and
$(y,z)\in\tau$. Now $(x,y)\in\alpha$ and $(y,z)\in\alpha$ imply that
$(x,z)\in\alpha$. Suppose that $(x,z)\notin\tau$, whence $(x,z)\in\alpha
\cap\alpha^{-1}$ and $[z]_{\gamma\cap\gamma^{-1}}<_{R}[x]_{\gamma\cap
\gamma^{-1}}$ follow. The relations $(y,z)\in\alpha$ and $(z,x)\in\alpha$
imply that $(y,x)\in\alpha$ and hence $(x,y)\in\alpha\cap\alpha^{-1}$.
Similarly, $(z,x)\in\alpha$ and $(x,y)\in\alpha$ imply that $(y,z)\in
\alpha\cap\alpha^{-1}$. In view of $(x,y)\in\tau$ and $(y,z)\in\tau$ we have
$[x]_{\gamma\cap\gamma^{-1}}\leq_{R}[y]_{\gamma\cap\gamma^{-1}}$ and
$[y]_{\gamma\cap\gamma^{-1}}\leq_{R}[z]_{\gamma\cap\gamma^{-1}}$, whence we
obtain that $[x]_{\gamma\cap\gamma^{-1}}\leq_{R}[z]_{\gamma\cap\gamma^{-1}}$,
a contradiction.

\noindent In order to prove that $\tau$ is a half-space we can use part (3) of
Proposition 2.2. Take $x,y,z\in A$\ such that $(x,y)\notin\tau$,
$(y,x)\notin\tau$\ and $(x,z)\in\tau$, $z\neq x$. Now $(x,y)\notin\tau$
implies that either $(x,y)\notin\alpha$ or $(x,y)\in\alpha\cap\alpha^{-1}$
with $[y]_{\gamma\cap\gamma^{-1}}<_{R}[x]_{\gamma\cap\gamma^{-1}}$. Similarly,
$(y,x)\notin\tau$ implies that either $(y,x)\notin\alpha$ or $(y,x)\in
\alpha\cap\alpha^{-1}$ with $[x]_{\gamma\cap\gamma^{-1}}<_{R}[y]_{\gamma
\cap\gamma^{-1}}$. It is easy to check that the only possibility to have
$(x,y)\notin\tau$ and $(y,x)\notin\tau$\ at the same time is the case when
$(x,y)\notin\alpha$ and $(y,x)\notin\alpha$. Since $\alpha$ is a half-space,
$(x,y)\notin\alpha$, $(y,x)\notin\alpha$ and $(x,z)\in\alpha$, $z\neq x$ imply
that $(y,z)\in\alpha$. Suppose that $(y,z)\in\alpha\cap\alpha^{-1}$, then
$(x,z)\in\alpha$ and the transitivity of $\alpha$ imply that $(x,y)\in\alpha$,
a contradiction. Thus we have $(y,z)\notin\alpha\cap\alpha^{-1}$, whence
$(y,z)\in\tau$ follows.$\square$

\bigskip

\noindent\textbf{2.5.Proposition.}\textit{\ Let the partial order }$\alpha
$\textit{\ be a half-space on }$A$\textit{. If }$\lambda$\textit{\ is a linear
order on }$A$\textit{, then}
\[
\alpha\lbrack\lambda]=\alpha\cup(\lambda\setminus(\alpha\cup\alpha^{-1}))
\]
\textit{is a linear extension of }$\alpha$\textit{ on }$A$\textit{\ and
}$\alpha=\alpha\lbrack\lambda]\cap\alpha\lbrack\lambda^{-1}]$\textit{.}

\bigskip

\noindent\textbf{Proof.} In order to see the transitivity of $\alpha
\lbrack\lambda]$ take $(x,y)\in\alpha\lbrack\lambda]$ and $(y,z)\in
\alpha\lbrack\lambda]$ with $x\neq y\neq z$. Clearly, $(x,y)\in\alpha$ and
$(y,z)\in\alpha$ imply $(x,z)\in\alpha$. If $(x,y)\in\alpha$ and
$(y,z)\in\lambda\setminus(\alpha\cup\alpha^{-1})$, then $(y,z)\notin\alpha$,
$(z,y)\notin\alpha$ and $(x,y)\in\alpha$, $x\neq y$, whence $(x,z)\in\alpha$
can be derived by part (4) of Proposition 2.2. Similarly, $(x,y)\in
\lambda\setminus(\alpha\cup\alpha^{-1})$ and $(y,z)\in\alpha$ imply
$(x,z)\in\alpha$ by part (3) of Proposition 2.2. If we have $(x,y)\in
\lambda\setminus(\alpha\cup\alpha^{-1})$ and $(y,z)\in\lambda\setminus
(\alpha\cup\alpha^{-1})$, then $(x,y)\in\lambda$ and $(y,z)\in\lambda$ imply
$(x,z)\in\lambda$. Since $(x,y)\notin\alpha\cup\alpha^{-1}$ and $(y,z)\notin
\alpha\cup\alpha^{-1}$ imply that $(x,y)\in\beta\cap\beta^{-1}$ and
$(y,z)\in\beta\cap\beta^{-1}$ (here $\beta$ is the complementary half-space of
$\alpha$), the transitivity of $\beta\cap\beta^{-1}$ gives that $(x,z)\in
\beta\cap\beta^{-1}$, i.e. that $(x,z)\notin\alpha\cup\alpha^{-1}$. It follows
that $(x,z)\in\lambda\setminus(\alpha\cup\alpha^{-1})$.

\noindent Suppose that $(x,y)\in\alpha\lbrack\lambda]$ and $(y,x)\in
\alpha\lbrack\lambda]$, then $(x,y)\in\alpha$ and $(y,x)\in\lambda
\setminus(\alpha\cup\alpha^{-1})$ is impossible. Similarly, $(x,y)\in
\lambda\setminus(\alpha\cup\alpha^{-1})$ and $(y,x)\in\alpha$ is also
impossible. Thus we have either $(x,y)\in\alpha$, $(y,x)\in\alpha$ or
$(x,y)\in\lambda\setminus(\alpha\cup\alpha^{-1})$, $(y,x)\in\lambda
\setminus(\alpha\cup\alpha^{-1})$, in both cases $x=y$ follows by the
antisymmetric properties of $\alpha$ and $\lambda$, respectively.

\noindent Suppose that $(x,y)\notin\alpha$ and $(y,x)\notin\alpha$, then
$(x,y)\notin\alpha\cup\alpha^{-1}$. Now $(x,y)\in\lambda$ implies
$(x,y)\in\lambda\setminus(\alpha\cup\alpha^{-1})$ and $(y,x)\in\lambda$
implies $(y,x)\in\lambda\setminus(\alpha\cup\alpha^{-1})$. We proved that
$\alpha\lbrack\lambda]$ is a linear order.

\noindent Using $\alpha\cap(\lambda\setminus(\alpha\cup\alpha^{-1}%
))=\alpha\cap(\lambda^{-1}\setminus(\alpha\cup\alpha^{-1}))=\varnothing$ and
$\lambda\cap\lambda^{-1}=\Delta_{A}$, it is straightforward to see that
$\alpha=\alpha\lbrack\lambda]\cap\alpha\lbrack\lambda^{-1}]$\textit{.}%
$\square$

\bigskip

\noindent\textbf{2.6.Corollary.}\textit{\ If }$\alpha$\textit{\ is a
half-space quasiorder on }$A$\textit{, then the induced partial order is of
the form }$r_{\alpha}=R_{1}\cap R_{2}$\textit{\ for some linear orders }%
$R_{1}$\textit{\ and }$R_{2}$\textit{\ on }$A/(\alpha\cap\alpha^{-1}%
)$\textit{, i.e. }$r_{\alpha}$\textit{ has order dimension at most }%
$2$\textit{.}

\bigskip

\noindent\textbf{Proof.} The partial order $r_{\alpha}$ is a half-space on
$A/(\alpha\cap\alpha^{-1})$ by Proposition 2.3. If $R$ is an arbitrary linear
order on $A/(\alpha\cap\alpha^{-1})$, then $R_{1}=r_{\alpha}[R]$ and
$R_{2}=r_{\alpha}[R^{-1}]$ are linear orders on $A/(\alpha\cap\alpha^{-1}%
)$\ with $r_{\alpha}=r_{\alpha}[R]\cap r_{\alpha}[R^{-1}]$ by Proposition
2.5.$\square$

\bigskip

\noindent We remark that Corollary 2.6 does not characterize half-spaces
entirely. As already noted, any linear order $\lambda$ on $A$ is an example of
a half-space: $\lambda\updownarrow\lambda^{-1}$. Let $f:A\longrightarrow X$ be
a function, $X_{1}\subseteq X$ a subset, $R$ a linear order on $X$ and define
the following relations on $A$:%
\[
\ker_{X_{1}}(f)=\Delta_{A}\cup\left\{  (a,b)\in A\times A\mid f(a)=f(b)\in
X_{1}\right\}  ,
\]%
\[
f^{-1}(R)=\Delta_{A}\cup\left\{  (a,b)\in A\times A\mid f(a)<_{R}f(b)\right\}
.
\]
The following is a standard construction of a half-space using a linear order.

\bigskip

\noindent\textbf{2.7.Proposition.}\textit{\ Let }$(A,\gamma)$\textit{\ be a
quasiordered set, }$(X,\rho)$\textit{\ a partially ordered set and
}$f:A\longrightarrow X$\textit{\ a }$(\gamma,\rho)$\textit{\ quasiorder
preserving function: }$(x,y)\in\gamma\Longrightarrow(f(x),f(y))\in\rho
$\textit{\ for all }$x,y\in A$\textit{. If }$X_{1}\subseteq X$\textit{\ is a
subset, }$\rho\subseteq R$\textit{\ is a linear extension of }$\rho$\textit{
on }$X$\textit{\ and }$\gamma\cap\ker(f)\subseteq\ker_{X_{1}}(f)$\textit{,
then}
\[
\alpha=\ker_{X_{1}}(f)\cup f^{-1}(R)
\]
\textit{is a half-space extension of }$\gamma$\textit{\ and }$\alpha\cap
\alpha^{-1}=\ker_{X_{1}}(f)$\textit{.}

\noindent\textit{If }$X_{1}=\varnothing$\textit{, then }$\ker_{X_{1}%
}(f)=\Delta_{A}$\textit{\ and }$\ker_{X_{1}}(f)\cup f^{-1}(R)=f^{-1}%
(R)$\textit{\ is a partial order.}

\noindent\textit{If }$X_{1}=X$\textit{, then }$\ker_{X_{1}}(f)=\ker
(f)$\textit{\ (now }$\gamma\cap\ker(f)\subseteq\ker_{X_{1}}(f)$%
\textit{\ automatically satisfied) and }$\ker_{X_{1}}(f)\cup f^{-1}%
(R)=\ker(f)\cup f^{-1}(R)$\textit{\ is a half-space extension of }$\gamma
$.\textit{\ In particular, if }$\kappa:A\longrightarrow A/(\gamma\cap
\gamma^{-1})$\textit{\ is the canonical surjection and }$R$\textit{\ is a
linear extension of the induced partial order }$r_{\gamma}$\textit{\ on
}$A/(\gamma\cap\gamma^{-1})$\textit{, then} $\ker(\kappa)\cup\kappa
^{-1}(R)=(\gamma\cap\gamma^{-1})\cup\kappa^{-1}(R)$\textit{\ is a half-space
extension of }$\gamma$\textit{.}

\bigskip

\noindent\textbf{Proof.} The containment $\gamma\subseteq\ker_{X_{1}}(f)\cup
f^{-1}(R)$ is a consequence of $\rho\subseteq R$, $\gamma\cap\ker
(f)\subseteq\ker_{X_{1}}(f)$ and of the quasiorder preserving property of $f$.
It is easy to see that $\ker_{X_{1}}(f)\cup f^{-1}(R)$ and $\ker_{X\setminus
X_{1}}(f)\cup f^{-1}(R^{-1})$ are quasiorders on $A$. We have
\[
(\ker_{X_{1}}(f)\cup f^{-1}(R))\cup(\ker_{X\setminus X_{1}}(f)\cup
f^{-1}(R^{-1}))=A\times A
\]
and%
\[
(\ker_{X_{1}}(f)\cup f^{-1}(R))\cap(\ker_{X\setminus X_{1}}(f)\cup
f^{-1}(R^{-1}))=\Delta_{A},
\]
thus $\ker_{X_{1}}(f)\cup f^{-1}(R)\updownarrow\ker_{X\setminus X_{1}}(f)\cup
f^{-1}(R^{-1})$. $\alpha\cap\alpha^{-1}=\ker_{X_{1}}(f)$ is obvious. To
conclude the proof, it is enough to note that $\kappa$ is a $(\gamma
,r_{\gamma})$ quasiorder preserving function.$\square$

\bigskip

\noindent\textbf{2.8.Proposition.}\textit{\ Let }$(A,\gamma)$\textit{\ be a
quasiordered set, }$(X,\rho)$\textit{\ a partially ordered set and
}$f:A\longrightarrow X$\textit{\ a completely }$(\gamma,\rho)$%
\textit{\ quasiorder preserving function: }$(x,y)\in\gamma\Longleftrightarrow
(f(x),f(y))\in\rho$\textit{\ for all }$x,y\in A$\textit{. If }$X_{1}%
^{(i)}\subseteq X$\textit{, }$i\in I$\textit{\ is a collection of subsets,
}$\gamma\cap\ker(f)\subseteq\ker_{X_{1}^{(i)}}(f)$\textit{\ for all }$i\in
I$\textit{\ and }$\left\{  R_{i}\mid i\in I\mathit{\ }\right\}  $\textit{\ is
a set of linear extensions of }$\rho$\textit{\ with }$\underset{i\in I}{\cap
}R_{i}=\rho$\textit{, then}
\[
\underset{i\in I}{\bigcap}(\ker_{X_{1}^{(i)}}(f)\cup f^{-1}(R_{i}))=\gamma,
\]
\textit{where the half-spaces }$\ker_{X_{1}^{(i)}}(f)\cup f^{-1}(R_{i}%
)$\textit{, }$i\in I$\textit{\ are described in Proposition 2.7. In
particular, if }$\kappa:A\longrightarrow A/(\gamma\cap\gamma^{-1}%
)$\textit{\ is the canonical surjection and }$\left\{  R_{i}\mid i\in
I\mathit{\ }\right\}  $\textit{\ is a set of linear extensions of the induced
partial order }$r_{\gamma}$\textit{\ on }$A/(\gamma\cap\gamma^{-1}%
)$\textit{\ with }$\underset{i\in I}{\cap}R_{i}=r_{\gamma}$\textit{, then}
\[
\underset{i\in I}{\bigcap(}\ker(\kappa)\cup\kappa^{-1}(R_{i}))=\underset{i\in
I}{\bigcap}((\gamma\cap\gamma^{-1})\cup\kappa^{-1}(R_{i}))=\gamma.
\]
\textbf{Proof.} We only have to show that
\[
\underset{i\in I}{\bigcap}(\ker_{X_{1}^{(i)}}(f)\cup f^{-1}(R_{i}%
))\subseteq\gamma.
\]
In view of the definition of $\ker_{X_{1}^{(i)}}(f)\cup f^{-1}(R_{i})$, the
relation%
\[
(a,b)\in\underset{i\in I}{\bigcap}(\ker_{X_{1}^{(i)}}(f)\cup f^{-1}(R_{i}))
\]
ensures that $f(a)\leq_{R_{i}}f(b)$ for all $i\in I$. Now $\underset{i\in
I}{\cap}R_{i}=\rho$ implies $f(a)\leq_{\rho}f(b)$, whence we obtain
$(a,b)\in\gamma$. To conclude the proof, it is enough to note that $\kappa$ is
completely $(\gamma,r_{\gamma})$ quasiorder preserving.$\square$

\bigskip

\noindent The following is now a straightforward consequence.

\bigskip

\noindent\textbf{2.9.Theorem.}\textit{\ Any quasiorder on }$A$\textit{\ can be
obtained as an intersection of half-space quasiorders on }$A$\textit{.}

\bigskip

\noindent In terms of the classification of convexities by separation axioms
(van de Vel [9]) the above theorem means that the convexity on $\{(x,y)\in
A\times A\mid x\neq y\}$ whose convex sets are the strict quasiorders on $A$
is an $S_{3}$ convexity, i.e. convex sets can be always separated from outside
points by complementary half-spaces, as in the standard convexity of an
Euclidean space or, as Szpilrajn's theorem [7] shows, in the coarser convexity
of strict partial orders plus $\{(x,y)\in A\times A\mid x\neq y\}$. However,
it is not difficult to see that, unlike in Euclidean space, in quasiorder
convexity, or in the coarser partial order convexity, disjoint convex sets
cannot always be separated by complementary half-spaces. A counterexample with
respect to both the quasiorder and partial order convexities is provided, for
$A=\{1,2,3,4\}$, by the partial orders $\{(1,2),(3,4)\}$ and $\{(1,4),(3,2)\}$.

\noindent Theorem 2.9 enables us to define a \textit{half-space realizer} of a
quasiorder $\gamma\subseteq A\times A$ as a set $\{\alpha_{i}\mid i\in I\}$ of
half-spaces on $A$ with $\underset{i\in I}{\bigcap}\alpha_{i}=\gamma$. The
\textit{half-space dimension} hs$\dim(A,\gamma)$ of a quasiordered set
$(A,\gamma)$ is the minimum of the cardinalities of the half-space realizers
of $\gamma$. The close analogy between the half-space dimension and the usual
order dimension of a partially ordered set can be seen immediately. The
observation preceding Proposition 2.2 guarantees that
\[
\text{hs}\dim(B,\gamma\cap(B\times B))\leq\text{hs}\dim(A,\gamma)
\]
for any subset $B\subseteq A$. Since any linear order is a half-space, for a
partially ordered set $(A,\gamma)$\ we have hs$\dim(A,\gamma)\leq\dim
(A,\gamma)$, where $\dim$ denotes the order dimension. In general, here we can
not expect equality. The partial order of the four element Boolean lattice
$M_{2}$ is a half-space, thus hs$\dim(M_{2},\leq)=1$, while $\dim(M_{2}%
,\leq)=2$. The next inequality is also a straightforward consequence of
Proposition 2.8.

\bigskip

\noindent\textbf{2.10.Corollary.}\textit{\ For a quasiordered set }%
$(A,\gamma)$\textit{\ we have}
\[
\text{hs}\dim(A,\gamma)\leq\dim(A/(\gamma\cap\gamma^{-1}),r_{\gamma}).
\]
\noindent The following theorem gives a complete description of half-space quasiorders.

\bigskip

\noindent\textbf{2.11.Theorem.}\textit{\ If }$\alpha\subseteq A\times
A$\textit{\ is a relation, then the following are equivalent.}

\begin{enumerate}
\item $\alpha$\textit{\ is a half-space quasiorder on }$A$\textit{.}

\item \textit{There exists an equivalence relation }$\varepsilon$\textit{\ on
}$A$\textit{, a linear order }$R$\textit{\ on the factor set }$A/\varepsilon
$\textit{\ and a function }$t:A/\varepsilon\longrightarrow\{0,1\}$%
\textit{\ with }$t([a]_{\varepsilon})=0$\textit{\ where }$[a]_{\varepsilon
}=\{a\}$\textit{\ such that}
\[
\alpha\!=\!\Delta_{A}\!\cup\!\left\{  (a,b)\!\in\!A\times A\!\mid
\![a]_{\varepsilon}=[b]_{\varepsilon}\text{ and }t([a]_{\varepsilon
})=1\right\}  \!\cup\!\left\{  (a,b)\!\in\!A\times A\!\mid\![a]_{\varepsilon
}<_{R}[b]_{\varepsilon}\right\}  .
\]

\item \textit{There exist a set }$X$\textit{, a subset }$X_{1}\subseteq
X$\textit{, a linear order }$R$\textit{ on }$X$\textit{ and a function
}$f:A\longrightarrow X$\textit{ such that }$\alpha=\ker_{X_{1}}(f)\cup
f^{-1}(R)$\textit{.}

\item \textit{There exists an equivalence relation }$\varepsilon$\textit{\ on
}$A$\textit{ such that }$\alpha$\textit{ is either the full or the identity
relation on each }$\varepsilon$\textit{-equivalence class, and any irredundant
set of representatives of the }$\varepsilon$\textit{-equivalence classes is
linearly ordered by }$\alpha$\textit{.}
\end{enumerate}

\bigskip

\noindent\textbf{Proof.} $(1)\Longrightarrow(2)$: Let $\alpha\updownarrow
\beta$ be complementary half-spaces and take
\[
\varepsilon=(\alpha\cap\alpha^{-1})\cup(\beta\cap\beta^{-1}).
\]
Clearly, $\varepsilon$ is reflexive and symmetric. Assume that $(x,y)\in
\alpha\cap\alpha^{-1}$ and $(y,z)\in\beta\cap\beta^{-1}$. Since $\alpha
\cup\beta=A\times A$, we have either $(x,z)\in\alpha$ or $(x,z)\in\beta$. In
the first case $(y,x)\in\alpha$ implies that $(y,z)\in\alpha\cap\beta
=\Delta_{A}$. In the second case $(z,y)\in\beta$ implies that $(x,y)\in
\alpha\cap\beta=\Delta_{A}$. Thus $(x,y)\in\alpha\cap\alpha^{-1}$ and
$(y,z)\in\beta\cap\beta^{-1}$ imply $x=y$ or $y=z$. Similarly, $(x,y)\in
\beta\cap\beta^{-1}$ and $(y,z)\in\alpha\cap\alpha^{-1}$ also imply $x=y$ or
$y=z$. In view of the above observations, it is easy to see that $\varepsilon$
is transitive. We also have $[a]_{\varepsilon}=[a]_{\alpha\cap\alpha^{-1}}%
\cup\lbrack a]_{\beta\cap\beta^{-1}}$ and $[a]_{\alpha\cap\alpha^{-1}}=\{a\}$
or $[a]_{\beta\cap\beta^{-1}}=\{a\}$\ for all $a\in A$.

\noindent We claim that $(a,b)\in\alpha$ and $[a]_{\varepsilon}\neq\lbrack
b]_{\varepsilon}$\ imply that $(x,y)\in\alpha$ for all $x\in\lbrack
a]_{\varepsilon}$ and for all $y\in\lbrack b]_{\varepsilon}$. Suppose that
$(x,y)\notin\alpha$, then $(x,y)\in\beta$. In view of $(x,a),(y,b)\in
\varepsilon$ we have the following cases. (i) $(x,a),(b,y)\in\alpha$, whence
$(x,y)\in\alpha$ can be obtained, a contradiction. (ii) $(x,a)\in\alpha$ and
$(y,b)\in\beta$, whence $(x,b)\in\alpha\cap\beta=\Delta_{A}$ can be obtained
in contradiction with $[x]_{\varepsilon}=[a]_{\varepsilon}\neq\lbrack
b]_{\varepsilon}$. (iii) $(a,x)\in\beta$ and $(b,y)\in\alpha$, whence
$(a,y)\in\alpha\cap\beta=\Delta_{A}$ can be obtained in contradiction with
$[a]_{\varepsilon}\neq\lbrack b]_{\varepsilon}=[y]_{\varepsilon}$. (iv)
$(a,x),(y,b)\in\beta$, whence $(a,b)\in\alpha\cap\beta=\Delta_{A}$ can be
obtained in contradiction with $[a]_{\varepsilon}\neq\lbrack b]_{\varepsilon}%
$. Thus the claim is proved.

\noindent Using our claim it is straightforward to check that%
\[
R=\left\{  ([a]_{\varepsilon},[b]_{\varepsilon})\mid(a,b)\in\alpha\right\}
\]
is a linear order on $A/\varepsilon$. For $a\in A$ let
\[
t([a]_{\varepsilon})=\left\{
\begin{array}
[c]{c}%
1\text{ if }[a]_{\varepsilon}=[a]_{\alpha\cap\alpha^{-1}}\neq\{a\}\\
0\text{\ otherwise \ \ \ \ \ \ \ \ \ \ \ \ \ \ \ \ \ }%
\end{array}
\right.  .
\]
Clearly, $t$ is well defined, moreover $[a]_{\varepsilon}=\{a\}$ implies
$[a]_{\varepsilon}=[a]_{\alpha\cap\alpha^{-1}}=\{a\}$ and $t([a]_{\varepsilon
})=0$. If $t([a]_{\varepsilon})=1$, then $[a]_{\varepsilon}=[a]_{\alpha
\cap\alpha^{-1}}$ and $[a]_{\varepsilon}=[b]_{\varepsilon}$\ implies that
$(a,b)\in\alpha$. It follows that%
\[
\Delta_{A}\!\cup\!\left\{  (a,b)\!\in\!A\times A\!\mid\![a]_{\varepsilon
}=[b]_{\varepsilon}\text{ and }t([a]_{\varepsilon})=1\right\}  \!\cup
\!\left\{  (a,b)\!\in\!A\times A\!\mid\![a]_{\varepsilon}<_{R}[b]_{\varepsilon
}\right\}  \!\subseteq\!\alpha.
\]
If $[a]_{\varepsilon}=[b]_{\varepsilon}$, then $(a,b)\in\alpha$ and $a\neq b$
implies that $[a]_{\varepsilon}=[a]_{\alpha\cap\alpha^{-1}}\neq\{a\}$, whence%
\[
\alpha\!\subseteq\!\Delta_{A}\!\cup\!\left\{  (a,b)\!\in\!A\times
A\!\mid\![a]_{\varepsilon}=[b]_{\varepsilon}\text{ and }t([a]_{\varepsilon
})=1\right\}  \!\cup\!\left\{  (a,b)\!\in\!A\times A\!\mid\![a]_{\varepsilon
}<_{R}[b]_{\varepsilon}\right\}
\]
can be obtained.

\noindent$(2)\Longrightarrow(3)$: It is straightforward to see that
$\alpha=\ker_{X_{1}}(f)\cup f^{-1}(R)$, where $X=A/\varepsilon$,
$X_{1}=\left\{  [a]_{\varepsilon}\mid a\in A\text{ and }t([a]_{\varepsilon
})=1\right\}  $ and $f:A\longrightarrow X$ is the canonical surjection. Thus
any half-space quasiorder can be obtained by the standard construction of
Proposition 2.7.

\noindent$(3)\Longrightarrow(1)$: This implication is a part of Proposition 2.7.

\noindent$(2)\Longleftrightarrow(4)$: Condition (4) is simply a reformulation
of (2).$\square$

\bigskip

\noindent\textbf{2.12.Remark.}\textit{ The triple }$(f,X_{1}\subseteq
X,R)$\textit{ given in the }$(2)\Longrightarrow(3)$\textit{ part of the above
proof has the following universal property. If }$g:A\longrightarrow Y$\textit{
is a function, }$Y_{1}\subseteq Y$\textit{ is a subset and }$S$\textit{\ is a
linear order on }$Y$\textit{ such that}%
\[
\alpha=\ker_{Y_{1}}(g)\cup f^{-1}(S),
\]
\textit{then there exists a unique function }$h:X\longrightarrow Y$\textit{
with }$h\circ f=g$\textit{, moreover }$h(X_{1})\subseteq Y_{1}$\textit{,
}$g^{-1}(\{y\})$\textit{ is a one element set for all }$y\in h(X\setminus
X_{1})\cap Y_{1}$\textit{ and }$h$\textit{\ is }$(<_{R},<_{S})$\textit{ strict
order preserving }

\bigskip

\noindent In view of the above characterization of the half-space $\alpha$, an
equivalence class $[a]_{\varepsilon}$ is called a \textit{box} of $\alpha$,
such a box is called \textit{full} if $t([a]_{\varepsilon})=1$ and
\textit{empty} if $t([a]_{\varepsilon})=0$ (note that a one element box is
always empty). A subset $B\subseteq A$ is a box of the half-space $\alpha$,
iff there are no elements $b_{1},b_{2}\in B$ such that $(b_{1},b_{2})\in
\alpha$, $(b_{2},b_{1})\notin\alpha$ and $B$ is maximal with respect to this
property. A box is empty if $\alpha\cap(B\times B)=\Delta_{B}$ and full if
$\left\vert B\right\vert >1$ and $B\times B\subseteq\alpha$.

\noindent In certain situations it is also convenient to give a half-space as%
\[
\alpha=(B_{w},w\in W,\leq_{W},t),
\]
where the subsets $B_{w}\subseteq A$, $w\in W$ are the boxes of $\alpha$, the
linear order $\leq_{W}$ is given on the index set $W$ and $t(B_{w})=1$ or
$t(B_{w})=0$ shows that $B_{w}$ is full or empty. If $W$ is finite, then we
can write $W=\{1,2,...,n\}$ and $\alpha=(B_{1}<B_{2}<...<B_{n},t)$. If
$\alpha\updownarrow\beta$ is a complementary pair of half-spaces, then
$\alpha$ and $\beta$ have the same boxes, a full $\alpha$-box is an empty
$\beta$-box and a full $\beta$-box is an empty $\alpha$-box, moreover the
linear order of the boxes in $\alpha$ and $\beta$ are opposite to each other.
It is also clear, that $[a]_{\alpha\cap\alpha^{-1}}=\{a\}$ if
$[a]_{\varepsilon}$ is empty and $[a]_{\alpha\cap\alpha^{-1}}=[a]_{\varepsilon
}$ if $[a]_{\varepsilon}$ is full.

\noindent With reference to the terminology of interval decompositions and
lexicographic sums of partial orders and more general relations (see e.g.
[3,4,5,6]), it is clear from condition (4) of Theorem 2.11 that half-space
quasiorders are precisely the lexicographic relational sums of trivial and
full binary relations over a linear order, i.e. they are the binary relations
decomposable into intervals such that the restriction to each interval is a
trivial or full relation and the quotient is a linear order.

\bigskip

\noindent\textbf{2.13.Theorem.}\textit{\ If }$(A,\gamma)$\textit{\ is a
quasiordered set and }$\{\alpha_{i}\mid i\in I\}$\textit{\ is a half-space
realizer of }$\gamma$\textit{ with }$\left|  I\right|  \geq2$\textit{, then
there exists an }$I$\textit{-indexed family }$R_{i}$\textit{, }$i\in
I$\textit{ of linear extensions of the induced partial order }$r_{\gamma}%
$\textit{ on }$A/(\gamma\cap\gamma^{-1})$\ \textit{such that}%
\[
\underset{i\in I}{\bigcap}R_{i}=r_{\gamma}.
\]
\noindent\textbf{Proof.} By Proposition 2.4, for each $i\in I$ there exists a
half-space $\tau_{i}$\ on $A$ such that $\gamma\subseteq\tau_{i}%
\subseteq\alpha_{i}$\ and $\tau_{i}\cap\tau_{i}^{-1}=\gamma\cap\gamma^{-1}$.
Clearly, $\underset{i\in I}{\cap}\alpha_{i}=\gamma$ implies that
$\underset{i\in I}{\cap}\tau_{i}=\gamma$, whence
\[
\underset{i\in I}{\bigcap}r_{\tau_{i}}=r_{\gamma}%
\]
can be derived for the induced partial orders $r_{\tau_{i}}$, $i\in I$ on
$A/(\tau_{i}\cap\tau_{i}^{-1})=A/(\gamma\cap\gamma^{-1})$. Using the notation
$\pi_{i}=r_{\tau_{i}}$, Proposition 2.3 ensures that each partial order
$\pi_{i}$ is a half-space on $P=A/(\gamma\cap\gamma^{-1})$.

\noindent We claim, that
\[
\rho=\Delta_{P}\cup\left(  \left(  \underset{i\in I}{\cup}\pi_{i}\right)
^{-1}\setminus\left(  \underset{i\in I}{\cup}\pi_{i}\right)  \right)
\]
is partial order on $P$. The reflexive and antisymmetric properties of $\rho$
can be immediately seen. In order to prove the transitivity of $\rho$ consider
the pairs $(x,y)\in\rho$ and $(y,z)\in\rho$ with $x,y,z\in P$ being different.
We have $(y,x)\in\pi_{j}$, $(z,y)\in\pi_{k}$ for some $j,k\in I$ and
$(x,y)\notin\underset{i\in I}{\cup}\pi_{i}$, $(y,z)\notin\underset{i\in
I}{\cup}\pi_{i}$. If $(z,y)\in\pi_{j}$, then the transitivity of $\pi_{j}$
implies $(z,x)\in\pi_{j}$. If $(z,y)\notin\pi_{j}$, then $(y,z)\notin\pi_{j}$
and the half-space property of $\pi_{j}$ imply that $(z,x)\in\pi_{j}$ (see
part (3) of Proposition 2.2). It follows that $(x,z)\in\left(  \underset{i\in
I}{\cup}\pi_{i}\right)  ^{-1}$. Suppose that $(x,z)\in\underset{i\in I}{\cup
}\pi_{i}$, then $(x,z)\in\pi_{t}$ for some $t\in I$. If $(z,y)\in\pi_{t}$,
then the transitivity of $\pi_{t}$ gives that $(x,y)\in\pi_{t}$, a
contradiction. If $(z,y)\notin\pi_{t}$, then $(y,z)\notin\pi_{t}$ and the
half-space property of $\pi_{t}$ gives that $(x,y)\in\pi_{t}$ (see part (4) of
Proposition 2.2), an other contradiction. Thus $(x,z)\notin\underset{i\in
I}{\cup}\pi_{i}$, whence $(x,z)\in\rho$ follows.

\noindent Let $\sigma_{i}\subseteq P\times P$ denote the complementary
half-space of $\pi_{i}$ and consider the following equivalence relation:%
\[
\Theta=\underset{i\in I}{\cap}(\sigma_{i}\cap\sigma_{i}^{-1})
\]
on $P$. Since $\pi_{i}^{-1}\cap\sigma_{i}^{-1}=\Delta_{P}$ for all $i\in I$,
we have $\rho\cap\Theta=\Delta_{P}$ and hence $\rho^{-1}\cap\Theta=\Delta_{P}%
$. Now we prove the containments $\Theta\circ\rho\subseteq\rho$ and $\rho
\circ\Theta\subseteq\rho$. If $(x,y)\in\Theta$ and $(y,z)\in\rho$ for the
elements $x,y,z\in P$ with $x,y,z$ being different, then $(z,y)\in\pi_{j}$ for
some $j\in I$ and $(y,z)\notin\underset{i\in I}{\cup}\pi_{i}$. In view of
$(x,y)\in\sigma_{j}\cap\sigma_{j}^{-1}$, we have $(x,y)\notin\pi_{j}$ and
$(y,x)\notin\pi_{j}$. Using part (4) in Proposition 2.2, we obtain that
$(z,x)\in\pi_{j}$ and $(x,z)\in\left(  \underset{i\in I}{\cup}\pi_{i}\right)
^{-1}$. Suppose that $(x,z)\in\underset{i\in I}{\cup}\pi_{i}$, then
$(x,z)\in\pi_{k}$ follows for some $k\in I$. Since $(x,y)\in\sigma_{k}%
\cap\sigma_{k}^{-1}$ implies that $(x,y)\notin\pi_{k}$ and $(y,x)\notin\pi
_{k}$, the application of part (3) in Proposition 2.2 yields $(y,z)\in\pi_{k}%
$, a contradiction. Thus we have $(x,z)\notin\underset{i\in I}{\cup}\pi_{i}$,
whence $(x,z)\in\rho$ follows. A similar argument shows that $\rho\circ
\Theta\subseteq\rho$.

\noindent Fix a linear order $\mu$ on $P$, then $\mu\cap\Theta$ and $\mu
^{-1}\cap\Theta$ are partial orders. Using the above properties of $\rho$ and
$\Theta$, it is straightforward to see that $\rho\cup(\mu\cap\Theta)$ and
$\rho\cup(\mu^{-1}\cap\Theta)$ are also partial orders on $P$.

\noindent Let $\rho\cup(\mu\cap\Theta)\subseteq\lambda$ and $\rho\cup(\mu
^{-1}\cap\Theta)\subseteq\lambda^{\ast}$ be linear extensions on $P$ and fix
an index $i^{\ast}\in I$. In view of Proposition 2.5, we can consider the
linear orders $R_{i}=\pi_{i}[\lambda]$, $i\in I\setminus\{i^{\ast}\}$ and
$R_{i^{\ast}}=\pi_{i^{\ast}}[\lambda^{\ast}]$ on $P$ (note that $I\setminus
\{i^{\ast}\}$\ is not empty). Since $\pi_{i}\subseteq R_{i}$ for all $i\in I$,
the inclusion
\[
r_{\gamma}=\underset{i\in I}{\bigcap}\pi_{i}\subseteq\underset{i\in I}%
{\bigcap}R_{i}%
\]
is obvious. In order to prove the reverse containment let $(x,y)\notin
\underset{i\in I}{\bigcap}\pi_{i}$ for some $x,y\in P$. We have $(x,y)\notin
\pi_{j}$ for some $j\in I$. If $(y,x)\in\underset{i\in I}{\cup}\pi_{i}$, then
$(y,x)\in\pi_{k}\subseteq R_{k}$ and hence $(x,y)\notin R_{k}$\ for some $k\in
I$. If $(y,x)\notin\underset{i\in I}{\cup}\pi_{i}$, then we distinguish two cases.

\noindent First suppose that $(x,y)\in\underset{i\in I}{\cup}\pi_{i}$. Then
$(y,x)\in\rho\subseteq\lambda\cap\lambda^{\ast}$ and the relations
$(x,y)\notin\pi_{j}$, $(y,x)\notin\pi_{j}$ imply that $(y,x)\in\pi_{j}%
[\lambda]$ (or $(y,x)\in\pi_{i^{\ast}}[\lambda^{\ast}]$ if $j=i^{\ast}$),
whence $(x,y)\notin R_{j}$ follows.

\noindent Next suppose that $(x,y)\notin\underset{i\in I}{\cup}\pi_{i}$. Then
$(x,y)\in\Theta$ and the linearity of $\mu$ gives that we have either
$(y,x)\in\mu\cap\Theta$ or $(y,x)\in\mu^{-1}\cap\Theta$. If $(y,x)\in\mu
\cap\Theta\subseteq\lambda$, then $(y,x)\in\pi_{i}[\lambda]$ and hence
$(x,y)\notin\pi_{i}[\lambda]=R_{i}$ for all $i\in I\setminus\{i^{\ast}\}$. If
$(y,x)\in\mu^{-1}\cap\Theta\subseteq\lambda^{\ast}$, then $(y,x)\in
\pi_{i^{\ast}}[\lambda^{\ast}]$ and hence $(x,y)\notin\pi_{i^{\ast}}%
[\lambda^{\ast}]=R_{i^{\ast}}$.$\square$

\bigskip

\noindent\textbf{2.14.Remark.}\textit{ Another possibility to construct the
linear orders }$R_{i}$\textit{ in the above proof is the following. Fix a well
ordering }$<$\textit{ on }$I$\textit{ and for }$i\in I\setminus\{i^{\ast}%
\}$\textit{ let}%
\[
R_{i}=\pi_{i}\cup((\sigma_{i}\cap\sigma_{i}^{-1})\cap\Lambda)\cup(\Theta
\cap\mu),
\]%
\[
R_{i^{\ast}}=\pi_{i^{\ast}}\cup((\sigma_{i^{\ast}}\cap\sigma_{i^{\ast}}%
^{-1})\cap\Lambda)\cup(\Theta\cap\mu^{-1}),
\]
\textit{where }$\Lambda=\{(x,y)\mid(y,x)\in\pi_{k}$ and $(x,y)\in
\underset{i\in I,i<k}{\bigcap}(\sigma_{i}\cap\sigma_{i}^{-1})$ for some $k\in
I\}$\textit{.}

\bigskip

\noindent In view of Corollaries 2.6 and 2.10, the above Theorem 2.13 yields
the following.

\bigskip

\noindent\textbf{2.15.Theorem.}\textit{\ If }$(A,\gamma)$\textit{\ is a
quasiordered set and }hs$\dim(A,\gamma)=1$\textit{, then }$\gamma$\textit{\ is
a half-space and}%
\[
\dim(A/(\gamma\cap\gamma^{-1}),r_{\gamma})=1\text{\textit{ if }}%
\gamma\text{\textit{ has no empty box with more than one element,}}%
\]%
\[
\dim(A/(\gamma\cap\gamma^{-1}),r_{\gamma})=2\text{\textit{ if }}%
\gamma\text{\textit{ has an empty box with more than one element.}}%
\]
\textit{If }hs$\dim(A,\gamma)\geq2$\textit{, then we have}%
\[
\dim(A/(\gamma\cap\gamma^{-1}),r_{\gamma})=\text{hs}\dim(A,\gamma).
\]
\noindent\textbf{2.16.Theorem.}\textit{\ If }$(A,\gamma)$\textit{\ is a
partially ordered set and }hs$\dim(A,\gamma)=1$\textit{, then }$\gamma
$\textit{\ is a half-space and}%
\[
\dim(A,\gamma)=1\text{\textit{ if }}\gamma\text{\textit{ is a linear order,}}%
\]%
\[
\dim(A,\gamma)=2\text{\textit{ if }}\gamma\text{\textit{ is not a linear
order.}}%
\]
\textit{If }hs$\dim(A,\gamma)\geq2$\textit{, then we have}%
\[
\dim(A,\gamma)=\text{hs}\dim(A,\gamma).
\]

\bigskip

\noindent3. DIRECT PRODUCT\ IRREDUCIBILITY\ OF\ HALF-SPACE QUASIORDERS

\bigskip

\noindent If $(A_{i},\gamma_{i})$, $i\in I$ is a family of quasiordered sets,
then
\[
\underset{i\in I}{\prod}\gamma_{i}=\{(\underline{a},\underline{b}%
)\mid\underline{a},\underline{b}\in\underset{i\in I}{\prod}A_{i}\text{ and
}(\underline{a}(i),\underline{b}(i))\in\gamma_{i}\text{ for all }i\in I\}
\]
is a quasiorder on the product set $\underset{i\in I}{\prod}A_{i}$ (here
$\underline{a}$ and $\underline{b}$ are functions $I\longrightarrow
\underset{i\in I}{\bigcup}A_{i}$\ such that $\underline{a}(i),\underline
{b}(i)\in A_{i}$ for all $i\in I$). We call $(\underset{i\in I}{\prod}%
A_{i},\underset{i\in I}{\prod}\gamma_{i})$ the direct product of the above
family. The kernel of the natural surjection%
\[
\varphi:\underset{i\in I}{\prod}A_{i}\longrightarrow\underset{i\in I}{\prod
}A_{i}/(\gamma_{i}\cap\gamma_{i}^{-1})
\]
is $\underset{i\in I}{\prod}(\gamma_{i}\cap\gamma_{i}^{-1})$, whence we obtain
a natural bijection%
\[
\left(  \underset{i\in I}{\prod}A_{i}\right)  /\left(  \underset{i\in I}%
{\prod}(\gamma_{i}\cap\gamma_{i}^{-1})\right)  \longrightarrow\underset{i\in
I}{\prod}A_{i}/(\gamma_{i}\cap\gamma_{i}^{-1}).
\]
It is easy to see that%
\[
(\underset{i\in I}{\prod}\gamma_{i})\cap(\underset{i\in I}{\prod}\gamma
_{i})^{-1}=\underset{i\in I}{\prod}(\gamma_{i}\cap\gamma_{i}^{-1})\text{ and
}r=\underset{i\in I}{\prod}r_{\gamma_{i}},
\]
where $r$ is the partial order on $\underset{i\in I}{\prod}A_{i}/(\gamma
_{i}\cap\gamma_{i}^{-1})$\ induced by the quasiorder $\underset{i\in I}{\prod
}\gamma_{i}$.

\noindent The product of non-trivial partial orders is never a linear order.
In contrast, the product of two half-spaces can be a half-space again: the
four element Boolean lattice $M_{2}$ is a product of two-element chains. We
show that this is the only possibility to get a non-trivial half-space as a
product of quasiorders.

\bigskip

\noindent\textbf{3.1.Lemma.}\textit{\ Let }$(A_{i},\gamma_{i})$\textit{,
}$i\in I$\textit{\ be a family of quasiordered sets and let }$j,k\in
I$\textit{, }$j\neq k$\textit{\ be indices such that }$a_{j}\neq c_{j}%
$\textit{, }$(a_{j},c_{j})\in\gamma_{j}$\textit{, }$(a_{j},b_{j})\notin
\gamma_{j}$\textit{\ for some }$a_{j},b_{j},c_{j}\in A_{j}$\textit{\ and
}$\gamma_{k}\neq A_{k}\times A_{k}$\textit{\ with }$\left|  A_{k}\right|
>1$\textit{. Then }$\underset{i\in I}{\prod}\gamma_{i}$\textit{\ is not a
half-space on }$\underset{i\in I}{\prod}A_{i}$\textit{.}

\bigskip

\noindent\textbf{Proof.} Let $\underline{u}\in\underset{i\in I}{\prod}A_{i}$
be an arbitrary element and $x_{k},y_{k}\in A_{k}$ such that $(x_{k}%
,y_{k})\notin\gamma_{k}$. Define $\underline{a},\underline{b},\underline{c}%
\in\underset{i\in I}{\prod}A_{i}$ as follows: for an index $i\in I$ let
\[
\underline{a}(i)=\left\{
\begin{array}
[c]{c}%
a_{j}\text{ if }i=j\text{ \ \ \ \ \ \ \ \ \ \ \ \ }\\
y_{k}\text{ if }i=k\text{ \ \ \ \ \ \ \ \ \ \ \ \ }\\
\underline{u}(i)\text{ if }i\in I\setminus\{j,k\}
\end{array}
\right.  ,\text{ \ }\underline{b}(i)=\left\{
\begin{array}
[c]{c}%
b_{j}\text{ if }i=j\text{ \ \ \ \ \ \ \ \ \ \ \ \ }\\
x_{k}\text{ if }i=k\text{ \ \ \ \ \ \ \ \ \ \ \ \ }\\
\underline{u}(i)\text{ if }i\in I\setminus\{j,k\}
\end{array}
\right.  ,
\]%
\[
\underline{c}(i)=\left\{
\begin{array}
[c]{c}%
c_{j}\text{ if }i=j\text{ \ \ \ \ \ \ \ \ \ \ \ \ }\\
y_{k}\text{ if }i=k\text{ \ \ \ \ \ \ \ \ \ \ \ \ }\\
\underline{u}(i)\text{ if }i\in I\setminus\{j,k\}
\end{array}
\right.  .
\]
Clearly, $(a_{j},b_{j})\notin\gamma_{j}$ implies $(\underline{a},\underline
{b})\notin\underset{i\in I}{\prod}\gamma_{i}$ and $(x_{k},y_{k})\notin
\gamma_{k}$ implies $(\underline{b},\underline{a})\notin\underset{i\in
I}{\prod}\gamma_{i}$. Since $(\underline{a},\underline{c})\in\underset{i\in
I}{\prod}\gamma_{i}$ and $(x_{k},y_{k})\notin\gamma_{k}$ implies
$(\underline{b},\underline{c})\notin\underset{i\in I}{\prod}\gamma_{i}$, we
can use part (3) in Proposition 2.2 to see that $\underset{i\in I}{\prod
}\gamma_{i}$ is not a half-space (we note that $\underline{c}\neq\underline
{a}$ is an immediate consequence of $a_{j}\neq c_{j}$).$\square$

\bigskip

\noindent\textbf{3.2.Lemma.}\textit{\ If }$(A,\gamma)$\textit{\ is a
quasiordered set such that there are no elements }$a,b,c\in A$\textit{\ with
}$a\neq c$\textit{, }$(a,c)\in\gamma$\textit{\ and }$(a,b)\notin\gamma
$\textit{, then }$\gamma\in\{\Delta_{A},A\times A\}$\textit{\ or }%
$\gamma=(B_{1}<B_{2})$\textit{\ is a half-space with a full lower box }$B_{1}%
$\textit{\ (or }$\left|  B_{1}\right|  =1$\textit{) and an empty upper box
}$B_{2}$\textit{.}

\bigskip

\noindent\textbf{Proof.} If $\gamma\notin\{\Delta_{A},A\times A\}$ satisfies
the above conditions, then for each $a\in A$ we have either $(a,x)\in\gamma$
for all $x\in A$ or $(a,y)\notin\gamma$ for all $y\in A$. Take
\[
B_{1}=\{a\in A\mid(a,x)\in\gamma\text{ for all }x\in A\}\text{ and }%
B_{2}=\{a\in A\mid(a,y)\notin\gamma\text{ for all }y\in A\},
\]
then $B_{1}\cup B_{2}=A$, $B_{1}\cap B_{2}=\varnothing$ and $\gamma
=B_{1}\times A=(B_{1}\times B_{1})\cup(B_{1}\times B_{2})$ is a half-space,
with a full lower box $B_{1}$ (or $\left|  B_{1}\right|  =1$) and an empty
upper box $B_{2}$. Thus we can write $\gamma=(B_{1}<B_{2})$.$\square$

\bigskip

\noindent\textbf{3.3.Lemma.}\textit{\ Let }$\gamma_{i}=(B_{i1}<B_{i2}%
)$\textit{, }$1\leq i\leq2$\textit{\ be half-spaces on }$A_{i}$\textit{\ with
full lower boxes }$B_{i1}$\textit{\ (or }$\left|  B_{i1}\right|  =1$\textit{)
and empty upper boxes }$B_{i2}$\textit{. Then we have the following.}

\begin{enumerate}
\item $\Delta_{A_{1}\times A_{2}}\neq\gamma_{1}\times\gamma_{2}\neq
(A_{1}\times A_{2})\times(A_{1}\times A_{2})$\textit{\ and take }%
$\underline{a}=(a_{12},a_{21})$\textit{, }$\underline{b}=(a_{11},a_{21}%
)$\textit{, }$\underline{c}=(a_{12},a_{22})$\textit{, where }$a_{ij}\in
B_{ij}$\textit{, }$i,j\in\{1,2\}$\textit{\ are arbitrary elements. Then
}$\underline{a}\neq\underline{c}$\textit{, }$(\underline{a},\underline{c}%
)\in\gamma_{1}\times\gamma_{2}$\textit{\ and }$(\underline{a},\underline
{b})\notin\gamma_{1}\times\gamma_{2}$\textit{.}

\item $\gamma_{1}\times\gamma_{2}$\textit{\ is a half-space if and only if
}$\left|  B_{ij}\right|  =1$\textit{\ for all }$i,j\in\{1,2\}$\textit{.}
\end{enumerate}

\bigskip

\noindent\textbf{Proof.}$(1)$: Obvious.

\noindent$(2)$: If $\left|  B_{ij}\right|  =1$ for all $i,j\in\{1,2\}$, then
it is clear that $A_{1}\times A_{2}$ is a four element set and $\gamma
_{1}\times\gamma_{2}$ is a partial order relation on $A_{1}\times A_{2}%
$\ providing a lattice isomorphic to $M_{2}$, which is a half-space as we have
already noted.

\noindent Suppose now, that $\left|  B_{11}\right|  >1$ and take $a^{\prime
},a^{\prime\prime}\in B_{11}$ such that $a^{\prime}\neq a^{\prime\prime}$. Let
$\underline{z}=(a^{\prime},b)$, $\underline{x}=(a^{\prime\prime},b)$ and
$\underline{y}=(a,c)$, where $a\in B_{12}$, $b\in B_{22}$, $c\in B_{21}$ are
arbitrary elements. Since $(\underline{x},\underline{y})\notin\gamma_{1}%
\times\gamma_{2}$, $(\underline{y},\underline{x})\notin\gamma_{1}\times
\gamma_{2}$ and $(\underline{x},\underline{z})\in\gamma_{1}\times\gamma_{2}$,
$(\underline{y},\underline{z})\notin\gamma_{1}\times\gamma_{2}$, we can apply
part (3) in Proposition 2.2 to derive that $\gamma_{1}\times\gamma_{2}$ is not
a half-space.

\noindent If $\left|  B_{12}\right|  >1$ then take $a^{\prime},a^{\prime
\prime}\in B_{12}$ such that $a^{\prime}\neq a^{\prime\prime}$. Let
$\underline{z}=(a^{\prime},b)$, $\underline{x}=(a^{\prime},c)$ and
$\underline{y}=(a^{\prime\prime},c)$, where $b\in B_{22}$, $c\in B_{21}$ are
arbitrary elements. Since $(\underline{x},\underline{y})\notin\gamma_{1}%
\times\gamma_{2}$, $(\underline{y},\underline{x})\notin\gamma_{1}\times
\gamma_{2}$ and $(\underline{x},\underline{z})\in\gamma_{1}\times\gamma_{2}$,
$(\underline{y},\underline{z})\notin\gamma_{1}\times\gamma_{2}$, we can apply
part (3) in Proposition 2.2 to derive that $\gamma_{1}\times\gamma_{2}$ is not
a half-space.

\noindent The cases $\left|  B_{21}\right|  >1$ and $\left|  B_{22}\right|
>1$ can be treated analogously.$\square$

\bigskip

\noindent\textbf{3.4.Theorem.}\textit{\ If }$(A_{i},\gamma_{i})$\textit{,
}$i\in I$\textit{\ is a family of non-trivial quasiordered sets (i.e. }%
$\Delta_{A_{i}}\neq\gamma_{i}\neq A_{i}\times A_{i}$\textit{\ for all }$i\in
I$\textit{), then the following are equivalent.}

\begin{enumerate}
\item $\underset{i\in I}{\prod}\gamma_{i}$\textit{\ is a half-space on
}$\underset{i\in I}{\prod}A_{i}$\textit{.}

\item \textit{Either }$I=\{1\}$\textit{\ and }$\gamma_{1}$\textit{\ is a
half-space or }$I=\{1,2\}$\textit{\ and }$(A_{1},\gamma_{1})$\textit{,
}$(A_{2},\gamma_{2})$\textit{\ are two-element chains.}
\end{enumerate}

\bigskip

\noindent\textbf{Proof.}

\noindent$(2)\Longrightarrow(1)$: It is an immediate consequence of part (2)
in Lemma 3.3.

\noindent$(1)\Longrightarrow(2)$: It is enough to deal with the case $\left|
I\right|  \geq2$. Using Lemma 3.1, we obtain that there is no $j\in I$ such
that $a_{j}\neq c_{j}$, $(a_{j},c_{j})\in\gamma_{j}$, $(a_{j},b_{j}%
)\notin\gamma_{j}$ for some $a_{j},b_{j},c_{j}\in A_{j}$. In view of Lemma
3.2, each $\gamma_{j}$ is a half-space on $A_{j}$\ of the form $\gamma
_{j}=(B_{j1}<B_{j2})$ with a full lower box $B_{j1}$ (or $\left|
B_{j1}\right|  =1$) and an empty upper box $B_{j2}$. If $\left|  I\right|
\geq3$, then we have different indices $i_{1},i_{2},i_{3}\in I$ and
\[
\underset{i\in I}{\prod}\gamma_{i}=(\gamma_{i_{1}}\times\gamma_{i_{2}}%
)\times\gamma_{i_{3}}\times(\underset{i\in I\setminus\{i_{1},i_{2},i_{3}%
\}}{\prod}\gamma_{i}),
\]
where $\gamma_{i_{1}}\times\gamma_{i_{2}}$ has the property described in part
(1) of Lemma 3.3. Since $\gamma_{i_{3}}\neq A_{i_{3}}\times A_{i_{3}}$ with
$\left|  A_{i_{3}}\right|  >1$, Lemma 3.1 ensures that our product is not a
half-space, a contradiction. Thus $\left|  I\right|  =2$ and part (2) in Lemma
3.3 gives that $(A_{1},\gamma_{1})$ and $(A_{2},\gamma_{2})$\ are two-element
chains (here we assumed $I=\{1,2\}$).$\square$

\bigskip

\noindent\textbf{3.5.Remark.}\textit{\ If }$\gamma_{j}=\Delta_{A_{j}}%
$\textit{\ and }$\left|  A_{j}\right|  >1$\textit{\ for some }$j\in
I$\textit{, then }$\underset{i\in I}{\prod}\gamma_{i}$\textit{\ is
disconnected, hence not a non-trivial half-space (because }$\underset{i\in
I}{\prod}\gamma_{i}=\Delta$\textit{\ would be the only possibility to get a
half-space). If }$\gamma_{j}=A_{j}\times A_{j}$\textit{\ for some }$j\in
I$\textit{, then }$\gamma_{j}$\textit{\ has no effect on wether the product
}$\underset{i\in I}{\prod}\gamma_{i}$\textit{\ is a half-space (in other words
}$\underset{i\in I}{\prod}\gamma_{i}$\textit{\ is a half-space if and only if
}$\underset{i\in I\setminus\{j\}}{\prod}\gamma_{i}$\textit{\ is a
half-space).}

\bigskip

\noindent\underline{Acknowledgement:} \textit{The initial version of this
paper was prepared while the first named author was at the Alfred Renyi
Institute of Mathematics, Hungarian Academy of Sciences.}

\bigskip

\noindent REFERENCES

\bigskip

\begin{enumerate}
\item Bonnet, R. , Pouzet, M. : \textit{Linear extensions of ordered sets,} in
Ordered Sets (I. Rival, ed.), Proceedings of the Nato Advanced Study Institute
Conference held in Banff, August 28-September 12, 1981, D. Reidel Publishing
Co., Dordrecht-Boston (1982), 125-170.

\item Dushnik, B., Miller, E.W. : \textit{Partially ordered sets,} Am. J.
Math. 63, 600-610 (1941).

\item Foldes, S. : \textit{On intervals in relational structures,} Zeitschrift
Math. Logik und Grundlagen Math. 26 (1980), 97-101.

\item Foldes, S. , Radeleczki, S. : \textit{On interval decomposition
lattices,} Discussiones Mathematicae, General Algebra and Applications 24
(2004), 95-114.

\item Hausdorff, F. : \textit{Grundz\"{u}ge einer Theorie der geordneten
Mengen,} Math. Ann. 65 (1908), no. 4, 435--505.

\item K\"{o}rtesi, P. , Radeleczki, S. , Szil\'{a}gyi, Sz.
:\textit{\ Congruences and isotone maps on partially ordered sets,} Math.
Pannonica 16/1 (2005), 39-55.

\item Szpilrajn, E. : \textit{Sur l'extension de l'ordre partiel,} Fund. Math.
16 (1930), 386-389.

\item Trotter, W. T.: \textit{Combinatorics and Partially Ordered Sets,
Dimension Theory,} The Johns Hopkins University Press, Baltimore-London, 1992.

\item van de Vel, M. L. J. : \textit{Theory of Convex Structures,}
North-Holland Mathematical Library, 50. North-Holland Publishing Co.,
Amsterdam, 1993.
\end{enumerate}

\end{document}